\documentclass[12pt]{amsart}
\usepackage{amsmath}
\usepackage{graphicx}
\usepackage{amsfonts}
\usepackage{amstext}
\usepackage{amsbsy}
\usepackage{amsopn}
\usepackage{amsxtra}
\usepackage{upref}
\usepackage{amsthm}
\usepackage{amsmath}
\usepackage{amssymb}

\def\D{{\mathbb D}}

\theoremstyle{definition}

\theoremstyle{remark}

\numberwithin{equation}{section}
\newcommand{\Pf}{\frac{f''}{f'}}
\newcommand{\noi}{\noindent}
\newcommand{\sm}{\smallskip}
\newcommand{\me}{\medskip}
\newcommand{\bi}{\bigskip}
\newcommand{\pf}{\me\noi{\bf Proof:} \,}
\begin{document}

\title {On Schwarz-Christoffel Mappings}

\author{Martin Chuaqui and Christian Pommerenke}
\thanks{The first
author was partially supported by Fondecyt Grant  \#1110321.
\endgraf  {\sl Key words:} Schwarz-Christoffel mapping, pre-vertices, convex, concave, univalent mapping, Blaschke product.
\endgraf {\sl 2000 AMS Subject Classification}. Primary: 30C20, 30C35;\,
Secondary: 30C45.}

\address{Facultad de Matem\'aticas\\ Pontificia Universidad Cat\'olica de Chile\\
Casilla 306, Santiago 22, CHILE.} \email{mchuaqui@mat.puc.cl}
\address{Institut f\"{u}r Mathematik MA 8-1\\Technische Universit\"{a}t\\
    D-10623 Berlin, GERMANY.} \email{pommeren@math.tu-berlin.de}
\email{}

\begin{abstract} We extend previous work on Schwarz-Chrsitoffel mappings, including the special cases when
the image is a convex polygon or its complement. We center our analysis on the relationship between the pre-Schwrazian
of such mappings and Blaschke products. For arbitrary Schwarz-Christoffel mappings, we resolve an open question in \cite{ChDO2}
that relates the degrees of the associated Blaschke products with the number of convex and concave vertices of the polygon.
In addition, we obtain a sharp sufficient condition in terms of the exterior angles for the injectivity of a mapping given
by the Schwarz-Christoffel formula, and study the geometric interplay between the location of the zeros of
the Blaschke products and the separation of the pre-vertices.

\end{abstract}
\maketitle

\section{Introduction}
The purpose of this paper is to provide further information about Schwarz-Christoffel mappings that add to the results
obtained in \cite{ChDO1}, \cite{ChDO2}. We refer the reader to \cite{BPW} for related interesting work on concave functions.

Let $f$ be a Schwarz-Christoffel mapping of the unit disk $\D$ onto the interior of an $(n+1)$-gon. In other words, $f$ is a conformal map onto a domain in the extended complex plane whose boundary consists of finitely many line segments, rays or lines. In \cite{ChDO2}, it is shown that the pre-Schwarzian
of $f$ has the form
\begin{equation}\label{Pf}\Pf = \frac{2B_1/B_2}{1-zB_1/B_2}
\end{equation}
for some finite Blaschke products $B_1, B_2$ without common zeros, with respective degrees $d_1, d_2$ satisfying $d_1+d_2=n$. The polygon is convex if and only if $d_2=0$ (see also \cite{ChDO1}).
The representation for $f''/f'$ is obtained from the well-known formula for such mappings,
\begin{equation}\label{SC-equation}\Pf=-2\sum_{k=1}^{n+1}\frac{\beta_k}{z-z_k} \, , \end{equation}
where $z_k$ are the pre-vertices and $2\pi\beta_k$ are the exterior angles with $\sum_{k=1}^{n+1}\beta_k=1$. We remind the reader of the relation
$2\pi\beta_k=\pi-\alpha_k$ between an exterior angle and the interior angle $\alpha_k$ at a given vertex. The Schwarz-Christoffel formula remains valid for polygons with one vertex at infinity. In this case the angle between the sides at infinity is, by definition, the angle (with minus sign) between the relevant sides (of their continuations) at a finite point.

As a consequence of equations (\ref{Pf}) and (\ref{SC-equation}),
the pre-vertices are shown to be the roots of the equation
\begin{equation}\label{roots} zB_1(z)/B_2(z)=1 \, .  \end{equation}
It is interesting that (\ref{roots}) corresponds to a polynomial equation of degree $n+1$ for which all roots are simple and lie on $|z|=1$. This is a particular feature of the
pair of Blaschke products $B_1, B_2$ appearing from Schwarz-Christoffel mappings. Note that the topological degree of $zB_1/B_2$ on $\partial\D$ is $1+d_1-d_2$,
so that $zB_1/B_2$ must be traversing in the negative sense at many of the pre-vertices. In fact, as the proof of Theorem 2  shows, at a pre-vertex $z_k$, $zB_1/B_2$ is traversing $\partial\D$
in the positive or negative sense according to whether $f(z_k)$ is a convex or a concave vertex. In is also interesting to observe that when $d_2=0$, any solution of (\ref{Pf}) will result in a univalent mapping because $1+\textrm{Re}\{zf''/f'\}\geq 0$.
In this paper we answer the natural question of finding a geometric interpretation for the degree $d_2$, and show that
this integer coincides with the number of concave vertices of the polygon.

\sm

In Section 2, we also address the case of Schwarz-Christoffel mappings $f$ onto the exterior of an $(n+2)$-gon,
with the normalization $f(0)=\infty$. In \cite{ChDO2}, we showed that
the pre-Schwarzian of such a mapping is given by
\begin{equation}\label{zPf}
z\Pf= \frac{2}{z^2(B_1/B_2)-1} \, ,
\end{equation}
for finite Blaschke products $B_1, B_2$  without common zeros, with degrees $d_1, d_2$, respectively, for which $d_1+d_2=n$. The polygon is convex if and only if $d_2=0$ and, as before, we show in this paper that $d_2$ is equal to the number of concave vertices of the polygon.

\sm

Another issue we address in this paper is the question of when a solution of (\ref{SC-equation}), or equivalently, of
$$ f'(z)=\prod_{k=1}^{n+1}(z-z_k)^{-2\beta_k} \;\, , \;\; \sum_{k=1}^{n+1}\beta_k=1 \, ,$$
does indeed correspond to
a {\it univalent} mapping. In Theorem 4 below we obtain the sharp sufficient condition $\sum_{k=1}^{n+1}|\beta_k|\leq 2$
for univalence. The result is optimal in the sense that there are non-univalent solutions of (\ref{SC-equation}) for which
$\sum_{k=1}^{n+1}|\beta_k|$ differs from $2$ by an arbitrarily small amount.

\sm

In Section 3 we obtain results on the separation of the pre-vertices of convex or concave
Schwarz-Christoffel mappings, expressed in terms of the location of the zeros $a_1,\ldots,a_n$ of the Blaschke
product $B_1$ that appears in (\ref{Pf}) or (\ref{zPf}) (recall that, in this case, $d_2=0$). The results are sharp, and show, for example, that the
pre-vertices tend to be uniformly separated on $\partial \D$ when all $|a_k|$ are very small.
Finally, in Section 4 we derive some necessary conditions for the location of the zeros of the Blaschke
products $B_1, B_2$  in (\ref{Pf}) and (\ref{zPf}) for arbitrary polygonal mappings.

\section{Blaschke Products and Univalence}

In \cite{ChDO1} we revisit the classical theme of convex mappings. The starting point is observation
that such mappings correspond exactly to the solutions of
$$ \Pf=\frac{2h}{1-zh} \, ,$$
for some function $h$ analytic in $\D$ and bounded by 1. The image $f(\D)$ is the interior of a polygon if and only if $h$ is a finite Blaschke product. We can express $h$ in terms of $p=f''/f'$ as
$$ h=\frac{p}{2+zp} \, ,$$
and draw the following result.

\bi\noi{\bf Theorem 1:} {\sl Let $h$ be analytic in $\D$ with $|h(z)|\leq 1$ everywhere. Then there exists a sequence $\{B_n\}_{n\in\mathbb{N}}$ of finite Blaschke products
converging to $h$ locally uniformly in $\D$.}

\pf Let $f$ be the convex mapping corresponding to $h$ as above, and let $\Omega_n$ be a sequence of convex polygons converging to $f(\D)$ in the sense of
Carath\'eodory. Properly normalized Schwarz-Christoffel mappings $f_n$ of $\D$ onto $\Omega_n$ will converge locally uniformly to $f$. Each mapping $f_n$ satisfies
(\ref{Pf}) for a certain finite Blaschke product $B_1=B_{1,n}$ and $B_2=1$. The theorem now follows by expressing $B_{1,n}$ in terms of the pre-Schwarzian of $f_n$.

\bi

Next, we give an answer to an important issue left unresolved in \cite{ChDO2}, namely the connection between the degrees $d_1,d_2$ and the number of convex and
concave vertices of the polygon.

\bi\noi {\bf Theorem 2:} {\sl Let $f$ map $\D$ onto the interior of an $(n+1)$-gon, and let $B_1, B_2$ be the corresponding Blaschke products
in the representation (\ref{Pf}). Then $d_2$ is equal to the number of concave vertices, while $d_1+1$ is equal to the number of convex vertices. }

\pf Let
\begin{equation}\label{boundary}\varphi(t)=\arg\left\{e^{it}{\frac{B_1}{B_2}(e^{it})}\right\} \, ,\end{equation}
with a well-defined branch of the argument once its value has been assigned at one given vertex. In any case,
\begin{equation}\label{boundary-1}\varphi'(t)=1+e^{it}\left(\frac{B_1'}{B_1}(e^{it})-\frac{B_2'}{B_2}(e^{it})\right)=1+|B_1'(e^{it})|-|B_2'(e^{it})| \, .\end{equation}
On the other hand, we see from (\ref{Pf}) and (\ref{SC-equation}) that
$$\frac{B_1/B_2}{zB_1/B_2-1}=\sum_{k=1}^{n+1}\frac{\beta_k}{z-z_k}\, ,$$
hence
\begin{equation}\label{angle}\beta_k=\lim_{z\rightarrow z_k}(z-z_k)\frac{B_1/B_2}{zB_1/B_2-1}=\frac{B_1/B_2}{(zB_1/B_2)'}(z_k)=\frac{1}{\varphi'(t_k)} \, ,
\end{equation}
where we have written $z_k=e^{it_k}$. We say that $z_k$ is convex or concave according to whether the polygon is convex or concave at $f(z_k)$. We conclude that $\varphi'(t_k)$
is positive at convex pre-vertices and negative at concave pre-vertices. Furthermore, the points $z_k$ are the solutions of equation (\ref{roots}). Hence we see from (\ref{boundary}) that, for all $k$,
$$ \varphi(t_k) = 2\pi j_k \; , \; j_k\in\mathbb{Z}\quad , \quad   \varphi(t)\neq 2\pi j \; , \; t\neq t_k \, .$$
It follows now  that
\sm
\begin{equation}\label{variation}\hspace{1cm} \int_{t_k}^{t_{k+1}}\varphi'(t)dt=\left\{\begin{array}{rl} 2\pi \;\; , & {\rm when} \: z_k, z_{k+1} \textrm{ are both convex} \\ & \\
0 \;\;\; ,& {\rm when} \: z_k, z_{k+1} \textrm{ are one convex and one concave } \hspace{1,3cm}  \\  & \\ -2\pi \; , & {\rm when} \: z_k, z_{k+1} \textrm{ are both concave} \: .  \end{array} \right.
\end{equation}

\me
\noi
Let $a$ be the number of consecutive convex pre-vertices, $b$ the number of instances a vertex of one type is followed by one of the other type, and $c$ the number of consecutive concave pre-vertices.
Then $a+b+c=n+1$, and we see by (\ref{variation}) that
$$ \int_0^{2\pi}\varphi'(t)dt=2\pi(1+d_1-d_2)=2\pi(a-c) \, .$$
Hence we have
$$1+d_1+d_2=n+1=a+b+c\quad , \quad 1+d_1-d_2=a-c \, .$$ We conclude that  that
$$1+d_1=a+\frac b2 \quad , \quad d_2=c+\frac b2 \, .$$
To obtain the theorem, we claim that $c+(b/2)$ is equal to the number of concave vertices (or pre-vertices). To see this, let $z_k,\ldots,z_l$ be any maximal chain of consecutive concave pre-vertices. Hence $z_{k-1}$ and $z_{l+1}$ are convex pre-vertices.
The collection $z_k,\ldots,z_l$ of concave pre-vertices contributes with $l-k$ in the count of $c$, and with 2 in the count of $b$. Thus its contribution in the count of $c+(b/2)$ is exactly the number of points in the chain. This  proves  the claim, and completes the proof
of the theorem.

\bi
Similar results hold for mappings $f$ onto the exterior of an $(n+2)$-gon, having the important normalization $f(0)=\infty$. For such mappings we have that
$$\Pf=2\left(\sum_{k=1}^{n+2}\frac{\beta_k}{z-z_k}-\frac1z\right)\, , $$
where, as before,  $z_k$ are the pre-vertices and $2\pi\beta_k$ are the exterior angles satisfying $-1<\beta_k<1$ and $\sum_{k=1}^{n+2}\beta_k=1$. In \cite{ChDO2},  this was shown
to lead to
$$z\Pf=\frac{2}{z^2(B_1/B_2)-1}  \, ,$$
for Blaschke products $B_1,B_2$ of degree $d_1,d_2$ satisfying $d_1+d_2=n$. Again,  the case $d_2=0$ corresponds exactly to when the polygon is convex. The pre-vertices appear as the solutions of the
equation $z^2B_1=B_2$, yet no further information was provided in connection with the degrees of the Blaschke products. With a similar argument as in the proof of Theorem 2, one can show:

\bi\noi {\bf Theorem 3:} {\sl Let $f$ map $\D$ onto the exterior of an $(n+2)$-gon, and let $B_1, B_2$ be the corresponding Blaschke products
in the representation (\ref{zPf}). Then $d_2$ is equal to the number of concave vertices, while $d_1+2$ is equal to the number of convex vertices.}

\me
Next, we address the question of the univalence of solutions of the equation (\ref{SC-equation}).

\bi\noi {\bf Theorem 4:} {\sl Let $0 \leq t_1 < \cdots < t_{n+1} < 2\pi,\,\,z_k =
e^{it_k},\,\,\beta_k \in \mathbb{R}$,\, $k=1,\ldots,n+1$, and let
\begin{equation}\label{optimal}
\sum^{n+1}_{k=1}\,\beta_k =1\quad , \quad \sum^{n+1}_{k=1}\,|\beta_k| \leq 2\, .
\end{equation}
Then the function $f$ defined by
\begin{equation}\label{SC-formal}
f'(z) = a \prod^{n+1}_{k=1}\,(z-z_k)^{-2\beta_k}\; , \, a\in\mathbb{C}\,,\, a\neq 0
\end{equation}
is univalent in $\D$.}

\me
Observe that there exist polygons with $\sum|\beta_k|$ arbitrarily large for which $f$ remains univalent. For
example, one can consider a polygon inscribed between two disjoint logarithmic spirals. On the other hand, once
$\sum^{n+1}_{k=1}\,|\beta_k|$ is allowed to exceed 2, then the sum of exterior angles at concave vertices will
be larger than $\pi$ in absolute value, making thus possible for the image $f(\D)$ to intersect itself.

\pf Let $f$ be given as in (\ref{SC-formal}) and suppose that (\ref{optimal}) holds. Then $f$ is locally injective in $\D$, and we will show that $f$ is univalent there,
in fact, that it is close-to-convex. Among the various equivalent formulations of this geometric property (see, {\it e.g.}, \cite{D83}, p. 48), we will
show for $0\leq\theta_1<\theta_2<2\pi$, that
$$I=I(\theta_1,\theta_2)=\int_{\theta_1}^{\theta_2}{\rm Re}\left\{1+z\frac{f''}{f'}(z)\right\}d\theta > -\pi \; , \; z=re^{i\theta} \, .$$
To prove this observe that
$$ {\rm Re}\left\{1+z\frac{f''}{f'}(z)\right\}=\sum_{k}\beta_k{\rm Re}\left\{\frac{z_k+z}{z_k-z}\right\}=\sum_{k}\beta_k\frac{1-r^2}{|z_k-z|^2} \, .$$
In trying to obtain a lower bound for $I$ we can discard the terms with $\beta_k>0$. For the other terms, we have that
$$\int_{\theta_1}^{\theta_2}\frac{1-r^2}{|z_k-re^{i\theta}|^2}d\theta < 2\pi \, , $$
because of the properties of the Poisson kernel. Hence
$$I > 2\pi\sum_{\beta_k<0}\beta_k \geq  -\pi \, ,$$
as desired.

\bi\noi{\bf Example:} Avkhadiev and Wirths initiated the study of the so called concave
mappings, that is, univalent mappings of the disk $\D$ onto the complement of a convex set
(see \cite{AW1}, \cite{AW2}, \cite{AW3}). As an example of Theorem 4 we can consider a
convex polygon $P$ with $\infty\in P$ and the conformal mapping of $\D$ onto the complement
of $P$. Let $\pi\lambda\,,1\leq\lambda\leq 2$, be the angle of
$f(\D)$ at $\infty$. It follows from  \cite{AW2} that
$$
f'(z) = a(z-z_{n+1})^{-\lambda-1}\,\prod_{k=1}^{n} (z-z_k)^{\gamma_k}\,,
\,\,\,\sum_{k=1}^{n}\gamma_k = \lambda-1 \, ,
$$
that is, $\beta_k = \frac12\gamma_k$ for $k=1,\ldots, n$, and $\beta_{n+1} =
\frac{1}{2}(1+\lambda)$. Therefore we have that
$$
\sum_{k=1}^{n+1} |\beta_{k}| = \frac{1}{2}(\lambda+1) +
\frac{1}{2}(\lambda-1) = \lambda \in [1,2] \, .
$$

Next, we establish the following variant of Theorem 4.

\bi\noi{\bf Theorem 5:} {\sl Let $f$ be defined by (\ref{SC-formal}) with $\sum_{k=1}^{n+1}\beta_k=1$, and suppose that}

\begin{equation}\label{symmetry}
\frac{\textrm{Im}\{f(z)\} }{\textrm{Im}\{z\} } >0\;\, \text{for
}|z|\leq1\, ,\,\text{Im }\{z\}\neq0.
\end{equation}
{\sl
Let $\theta_{\pm}$ be the interior angles of the polygon
$f(\partial\D)$ at $f(\pm1)$. If
\begin{equation}\label{optimal-1}
\sum_{k=1}^{n+1} |\beta_{k}|\,
\leq\, 3 + \frac{1}{\pi} \max(\theta_{+}-\pi,0) + \frac{1}{\pi} \max
(\theta_{-}-\pi,0)
\end{equation}
then $f$ is univalent.}

\bi
The expression in the right hand side of \eqref{optimal-1} lies in $[3,5]
$, and it is easy
to see that any value in this range can be achieved. Therefore, Theorem 5 gives a better
result than Theorem 4 under the stronger assumption \eqref{symmetry}. The condition \eqref{symmetry}
implies, in particular, that $f(\D)$ is symmetric with respect to $\mathbb{R}$.

\pf By (\ref{symmetry}) the polygon $P=f(\partial\D)$ is symmetric with respect to $\mathbb{R}$. Hence $m:=(n+1)/2\in\mathbb{N}$.
We may assume that $z_1=1, z_2=-1$ in (\ref{SC-formal}). Then
$$
\beta_1 = \frac{1}{2} - \frac{\theta_{+}}{2\pi}\quad , \quad
\beta_2 = \frac{1}{2} - \frac{\theta_{-}}{2\pi}\, ,
$$
which satisfy $|\beta_1|, |\beta_2|\leq\frac12$.
It follows that
$$
\frac{1}{\pi}\max(\theta_{+}-\pi,0)=\max(-2\beta_1,0)=|\beta_1|-\beta_1 \, ,
$$
$$
\frac{1}{\pi}\max(\theta_{-}-\pi,0)=\max(-2\beta_2,0)=|\beta_2|-\beta_2 \, .
$$

Let $\varphi_\pm$ be the conformal mappings of $\D$ onto the semi-discs
$\{z\in\D : \text{Im }z \gtrless 0 \}$ such that $\varphi_\pm (1)=1
, \varphi_\pm (-1)=-1 $
and $\varphi_\pm (\pm i) = \pm i$ . Then
\begin{equation}\label{poligonos}
P_{\pm}=f(\varphi_{\pm}(\partial\D)) = f(\partial\D \cap\{\text{Im
}z \gtrless 0 \})\, \cup\, [f(-1),f(+1)]
\end{equation}
are the upper and lower parts of $P$ union $[f(-1),f(1)]$. We may also assume that $\beta_k\, , \, k=3,\ldots, m+1$,
belong to the vertices of $P$ that lie in $P_{+}$.

Let us consider the upper polygon $P_{+}$. The values $\gamma_k$ of $P_{+}$ corresponding to the $\beta_k$
are
$$\gamma_k=\left\{\begin{array}{ll} \displaystyle{\frac14+\frac12\beta_k}\geq 0 \; , \; k=1,2\\ \\ \beta_k \; , \; k=3,\ldots, m+1 \, ,\end{array} \right .
$$
(for which $\sum_{k=1}^{m+1}\gamma_k=1$.)
In light of the symmetry with respect to $\mathbb{R}$ we get
$$
2\sum_{k=1}^{m+1}|\gamma_k|=1+\beta_1+\beta_2+\sum_{k=3}^{n+1}|\beta_k|=1+\sum_{k=1}^{n+1}|\beta_k|-\left(|\beta_1|-\beta_1\right)
-\left(|\beta_2|-\beta_2\right) \, .$$
Using (\ref{optimal-1}) we conclude that
$$ 2\sum_{k=1}^{m+1}|\gamma_k| \leq 4 \, ,$$
and it follows from Theorem 4 that $f\circ\varphi_{+}$ is univalent in $\D$. The same holds for $f\circ\varphi_{-}$.

By (\ref{poligonos}) we have
$$ f(\D)=(f\circ\varphi_{+})(\D)\cup f((-1,1))\cup (f\circ\varphi_{-})(\D) \, ,$$
which are disjoint unions by (\ref{symmetry}). Hence $f$ is univalent in $\D$.

\section{Separation of Pre-Vertices}

Let $f$ be a Schwarz-Christoffel mapping taking $\D$ onto a convex $(n+1)$-gon. Recall that
$$ \frac{f''}{f'}(z)=\frac{2B(z)}{1-zB(z)} \, ,$$
where $B(z)$ is a Blaschke product of degree $n$. We write
$$B(z)=c\prod_{k=1}^n\frac{z-a_k}{1-\overline{a_k}z}  \, ,
$$
where $|c|=1$ and all $|a_k|<1$. The pre-vertices $z_1,\ldots, z_{n+1}$ correspond to the roots of the equation $zB(z)=1$,
and after rotating $f$, we may assume that $c=1$. Recall also, that in this case, any choice
of Blaschke product $B=B(z)$ will result in a univalent mapping $f$. The separation of consecutive pre-vertices
$z_k,z_{k+1}$ is to be understood as $\arg\{\bar{z}_kz_{k+1}\}\in (0,2\pi)$.

\me\noi{\bf Theorem 6:} {\sl
Suppose that $|a_k|\leq r <1$ for all $k$. Then

\me\noi
(i) The minimum separation in argument of consecutive pre-vertices is given by $2\theta$, where $\theta$ is the
unique root in $(0,\pi/2)$ of
the equation
\begin{equation}\label{min} \pi=\theta+2n\arctan\left\{\frac{1+r}{1-r}\,\tan\left(\frac{\theta}{2}\right)\right\} \, .
\end{equation}
The result is sharp. The optimal configuration occurs when $a_1=\cdots=a_n=rc$
for some root of the equation $c^{n+1}=-1$, and the lower bound is attained for the pre-vertices
$$ e^{i\theta}c\;\, ,\;e^{-i\theta}c \, .$$
The distance between any other pair of consecutive pre-vertices will be larger.

\me\noi
(ii) The maximum separation in argument of consecutive pre-vertices is given by $2\psi$, where $\psi$ is the
unique root in $(0,\pi)$ of
the equation
\begin{equation}\label{max} \pi=\psi+2n\arctan\left\{\frac{1-r}{1+r}\,\tan\left(\frac{\psi}{2}\right)\right\} \, .\end{equation}
The result is sharp. The optimal configuration occurs when $a_1=\cdots=a_n=rd$
for some root of the equation $d^{n+1}=(-1)^n$, and the upper bound is attained for the pre-vertices
$$ -e^{i\psi}d\;\, ,\;-e^{-i\psi}d \, .$$
The distance between any other pair of consecutive pre-vertices will be larger. }

\pf We need to estimate the distance between two consecutive roots $a=e^{i\alpha},\, b=e^{i\beta}$ of the equation $zB(z)=1$.
Because $zB(z)$ traces the boundary $\partial\D$ for $z\in\partial\D$ in a monotonic way, we must have from (\ref{boundary-1}) and (\ref{variation}) that
\begin{equation}\label{2pi} \int_{\alpha}^{\beta}\left(1+|B'(e^{it})|\right)dt=2\pi \, ,
\end{equation}
with
$$|B'(e^{it})|=\sum_{k=1}^n\frac{1-|a_k|^2}{|e^{it}-a_k|^2} \, . $$
(Equation (\ref{2pi}) shows that $0<\beta-\alpha<2\pi$.) We claim that for $\alpha, \beta$ fixed, the contribution of any single
summand
\begin{equation}\label{integral} \int_{\alpha}^{\beta}\frac{1-|a_k|^2}{|e^{it}-a_k|^2}dt
\end{equation}
will be maximal if $|a_k|=r$ and $a_k/|a_k|$ is the midpoint $c$
of the shorter arc joining $a$ and $b$. Let $r_k=|a_k|$ and write $a_k=r_ke^{it_k}$. Then
$$ \int_{\alpha}^{\beta}\frac{1-|a_k|^2}{|e^{it}-a_k|^2}dt=\int_{\alpha-t_k}^{\beta-t_k}\frac{1-r_k^2}{|1-r_ke^{it}|^2}dt \, .$$
For $r_k\leq r$ given, this integral is maximal when $1\in\partial\D$ is the midpoint of the shorter arc between
$e^{i(\alpha-t_k)}$ and $e^{i(\beta-t_k)}$. The integral is then equal to
$$ \int_{-\theta}^{\theta}\frac{1-r_k^2}{|1-r_ke^{it}|^2}dt =4\arctan\left\{\frac{1+r_k}{1-r_k}\tan\left(\frac{\theta}{2}\right)\right\} \, , $$
where $2\theta=\beta-\alpha$, and becomes maximal if $r_k=r$. In other words,
\begin{equation}\label{maximal} \int_{\alpha}^{\beta}\frac{1-|a_k|^2}{|e^{it}-a_k|^2}dt\, \leq \,
4\arctan\left\{\frac{1+r}{1-r}\tan\left(\frac{\theta}{2}\right)\right\} \, ,
\end{equation}
which proves our claim for the contribution of any single term, and therefore
the minimum separation between consecutive roots will
occur if this holds for all $k=1, \ldots,n$.
Equation (\ref{min}) follows. The analysis shows that for the extremal configuration, all $a_k=rc$ are equal and that
$e^{\pm i\theta}c$ are roots of the equation $zB(z)=1$. Because $B(c)=c$ then $cB(c)=c^{n+1}$, and since $zB(z)$
traces the arc between the two roots in symmetric fashion with respect to the midpoint, we conclude that $c^{n+1}=-1$.
This proves part (i).

\sm
For part (ii), we observe that, for $r_k=|a_k|$ fixed, (\ref{integral}) will be minimal provided $-1\in\partial\D$ is the midpoint of the shorter arc between $e^{i(\alpha-t_k)}$ and $e^{i(\beta-t_k)}$. The integral is then equal to
$$ \int_{\pi-\psi}^{\pi+\psi}\frac{1-r_k^2}{|1-r_ke^{it}|^2}dt =4\arctan\left\{\frac{1-r_k}{1+r_k}\tan\left(\frac{\psi}{2}\right)\right\} \, , $$
where $2\psi=\beta-\alpha$, and becomes minimal when $r_k=r$. Thus,
\begin{equation}\label{minimal} \int_{\alpha}^{\beta}\frac{1-|a_k|^2}{|e^{it}-a_k|^2}dt\, \geq \,
4\arctan\left\{\frac{1-r}{1+r}\tan\left(\frac{\theta}{2}\right)\right\} \, .
\end{equation}
Therefore, the maximum separation between
consecutive roots will occur if, for all $k=1,\ldots,n$, we have that $|a_k|=r$ and $a_k/r$ is equal to the midpoint of the longer arc between $a$ and $b$. From this, equation (\ref{max}) follows. As before, the
analysis of the extremal configuration gives $a_k=rd$ for all $k$. The points $-e^{\pm i\psi}d$ are roots of $zB(z)=1$,
which by symmetry as before implies that $d^{n+1}=(-1)^n$.

\bi
\noi{\bf Corollary 7:} {\sl Suppose that  $|a_k|\leq \epsilon$ for all $k$. Then the maximum separation $2\psi$
and minimum separation $2\theta$
between consecutive pre-vertices satisfy
\begin{equation}\label{concentration}
\frac{\pi}{1+(1+2\epsilon)n}+O(\epsilon^2)\,\leq\, \theta\,\leq\, \psi\, \leq\,
\frac{\pi}{1+(1-2\epsilon)n}+O(\epsilon^2)\; , \;\, \epsilon\rightarrow 0 \, .
\end{equation}}

\pf For fixed $x\in[0,\pi/2]$, let $F(\delta)=\arctan\left((1+\delta)\tan(x)\right)$. Then $F(0)=x$, $F'(0)=\sin(x)\cos(x)=\frac12\sin(2x)$ and
$F''(0)=-2\sin^3(x)\cos(x)$, hence
$$F(\delta)=x+\frac12\sin(2x)\delta+O(\delta^2)\; , \;\, \delta\rightarrow 0 \, .$$
Using that $(1+r)/(1-r)=1+2r+O(r^2), r\rightarrow 0$, and that $\sin(2x)\leq 2x$, we see from (\ref{min}) that the
minimum separation $\theta$ satisfies
$$
\pi\leq\theta+2n\left(\frac{\theta}{2}+\epsilon\sin(\theta)+O(\epsilon^2)\right)\leq
\left(1+(1+2\epsilon)n\right)\theta+O(\epsilon^2) \, .$$
This implies the lower bound in (\ref{concentration}).
A similar analysis applies to the maximum separation $\psi$, and the upper bound in (\ref{concentration}) obtains.

\bi
Suppose now that $f$ is a Schwarz-Christoffel mapping taking $\D$ onto the complement of a bounded convex $(n+2)$-gon, with the
normalization $f(0)=\infty$.  In this situation, we know that
$$ z\frac{f''}{f'}(z)=\frac{2}{z^2B(z)-1} \, ,$$
where $B(z)$ again is a Blaschke product of degree $n$. We write
$$ B(z)=c\prod_{k=1}^n\frac{z-a_k}{1-\overline{a_k}z}  \, ,$$
where $|c|=1$ and all $|a_k|<1$.
The pre-vertices $z_1,\ldots, z_{n+2}$ are now given by the roots
of the equation $z^2B(z)=1$, and after a rotation of $f$, we may assume that $c=1$. The following result is obtained
in a way similar to Theorem 6, and the proof will be omitted.

\me\noi{\bf Theorem 8:} {\sl
Suppose that $|a_k|\leq r <1$ for all $k$. Then

\me\noi
(i) The minimum separation in argument of consecutive pre-vertices is given by $2\theta$, where $\theta$ is the
unique root in $(0,\pi/2)$ of
the equation
$$ \pi=2\theta+2n\arctan\left\{\frac{1+r}{1-r}\tan\left(\frac{\theta}{2}\right)\right\} \, .$$
The result is sharp. The optimal configuration occurs when $\alpha_1=\cdots=\alpha_n=rc$
for some root of the equation $c^{n+2}=-1$, and the lower bound is attained for the pre-vertices
$$ e^{i\theta}c\;\, ,\;e^{-i\theta}c \, .$$
The distance between any other pair of consecutive pre-vertices will be larger.

\me\noi
(ii) The maximum separation in argument of consecutive pre-vertices is given by $2\psi$, where $\psi$ is the
unique root in $(0,\pi/2)$ of
the equation
$$ \pi=2\psi+2n\arctan\left\{\frac{1-r}{1+r}\tan\left(\frac{\psi}{2}\right)\right\} \, .$$
The result is sharp. The optimal configuration occurs when $\alpha_1=\cdots=\alpha_n=rd$
for some root of the equation $d^{n+2}=(-1)^{n+1}$, and the upper bound is attained for the pre-vertices
$$ -e^{i\psi}d\;\, ,\;-e^{-i\psi}d \, .$$
The distance between any other pair of consecutive pre-vertices will be larger. }

\bi A statement similar to Corollary 7 can be made in this case. If $|a_k|\leq \epsilon$ then
the maximum and minimum separation between pre-vertices satisfy
\begin{equation}\label{concentration-1}
\frac{\pi}{2+(1+2\epsilon)n}+O(\epsilon^2)\,\leq\, \theta\,\leq\, \psi\, \leq\,
\frac{\pi}{2+(1-2\epsilon)n}+O(\epsilon^2)\; , \;\, \epsilon\rightarrow 0 \, .
\end{equation}

\bi

We finish this section with some remarks on the separation of pre-vertices for arbitrary polygonal mappings. Suppose $f$ is a mapping
of the form given by (\ref{Pf}), where after rotation, we may assume expressions for $B_1, B_2$ given by
$$B_1(z)=\prod_{k=1}^{d_1}\frac{z-a_k}{1-\overline{a_k}z} \quad , \quad B_2(z)=\prod_{k=1}^{d_2}\frac{z-b_k}{1-\overline{b_k}z} \; .$$
Then
$$\varphi'(t)=1+|B_1'(e^{it})|-|B_2'(e^{it})|=1+\sum_{k=1}^{d_1}\frac{1-|a_k|^2}{|e^{it}-a_k|^2}
-\sum_{k=1}^{d_2}\frac{1-|b_k|^2}{|e^{it}-b_k|^2}\, .$$
Let $a=e^{i\alpha}, b=e^{i\beta}$ be consecutive {\it convex} pre-vertices, with separation $\beta-\alpha=2\delta$,
and let $r$ be the radius of the smallest centered subdisk that contains the zeros of $B_1, B_2$.
We deduce from (\ref{variation}) and the estimates (\ref{maximal}), (\ref{minimal}), that
\begin{equation}\label{sep-convex}
\end{equation}
$$ \delta+2d_1\arctan\left(x/\lambda\right)-2d_2\arctan\left(\lambda x\right)\, \leq \, \pi \,\leq \,
\delta+2d_1\arctan\left(\lambda x\right)-2d_2\arctan\left(x/\lambda\right) \, ,$$

\me
\noi
where $\lambda=(1+r)/(1-r)$ and $x=\tan\left(\delta/2\right)$. Thus, for example, with given $d_1, d_2$, a
relatively small separation $2\delta$ can only occur if $r$ is rather close to 1.
Because the univalence of $f$ is no longer guaranteed when $B_1, B_2$ are chosen arbitrarily,
is seems of interest to determine under which circumstances the inequalities (\ref{sep-convex}) remain sharp. We provide here
a simple example where one can show sharpness in the right-hand side of (\ref{sep-convex}) when $d_1=d_2=1$.

\me
\noi{\bf Example:} Consider the Blaschke products $B_1, B_2$ given by
$$ B_1(z)=\frac{z+r}{1+rz}\quad , \quad B_2(z)=\frac{z-r}{1-rz} \;\, , \; r\in(0,1) \, ,$$
and let $f$ be defined, up to an affine change, by
$$ \frac{f''}{f'}=\frac{2B_1/B_2}{1-zB_1/B_2}=\frac{2(z+r)(1-rz)}{(z-r)(1+rz)-z(z+r)(1-rz)} \, .$$
In analyzing the roots of $zB_1=B_2$, we observe that $z_3=1$ is one immediate solution. The other
solutions are the roots of
$$rz^2+(r^2+2r-1)z+r=0\, ,$$
which are given by
$$z_{1,2}=\frac{(1-2r-r^2)\pm\sqrt{-(1-r^2)(r^2+4r-1)}}{2r}\, .$$
For $r> r_0=\sqrt{5}-2=0,236\ldots$, the discriminant is negative and $|z_{1,2}|=1$, with $z_1=z_2(=-1)$ only for $r=1$.
In the partial fraction decomposition
$$\frac{f''}{f'}=-2\left(\frac{\beta_1}{z-z_1}+\frac{\beta_2}{z-z_2}+\frac{\beta_3}{z-z_3}\right) \, ,$$
we must have $\beta_1=\beta_2$ because of symmetry, while $\beta_1+\beta_2+\beta_3=1$ by equating
coefficients with the above representation for $f''/f'$. Recall equation (\ref{angle}) that relates
the exterior angles $2\pi\beta_k$ with the boundary function $\varphi(t)$. Here
$$\varphi'(t)=1+\frac{1-r^2}{|1+re^{it}|^2}-\frac{1-r^2}{|1-re^{it}|^2} \, ,$$
hence
$$\varphi'(0)=1+\frac{1-r}{1+r}-\frac{1+r}{1-r} \, .$$
One readily verifies that $\varphi'(0)\leq -2$ precisely when $r\geq r_1=(1+\sqrt{13})/(5+\sqrt{13})=0,535\ldots$,
in which case $\beta_3\in(-\frac12,0)$ and $\beta_1=\beta_2\in(\frac12,\frac34)$. Thus, for  $r\geq r_1$, we deduce
from Theorem 4 that $f$ is {\it univalent}, and $z_{1,2}$ are convex pre-vertices
while $z_3$ is a concave pre-vertex. The convex vertices $f(z_{1,2})$ are at infinity, and the image
$f(\D)$ corresponds to a half-plane minus a symmetric slit ending at the concave vertex when $r=r_1$, or a wedge when $r>r_1$.

\sm

Finally, to show sharpness in the right-hand side of (\ref{sep-convex}), observe that for $r\geq r_1$ the conjugate points
$z_1=\bar{z_2}$ have negative real part, and thus their separation $2\theta$ will correspond to the root $\theta\in(0,\pi/2)$
of the equation
$$\theta+2\arctan(\lambda x)-2\arctan(\lambda^{-1}x)=\pi \, .$$

\sm

For consecutive {\it concave} pre-vertices we deduce in similar fashion that
\begin{equation}\label{sep-concave}
\end{equation}
$$\delta+2d_1\arctan\left(x/\lambda\right)-2d_2\arctan\left(\lambda x\right)\, \leq \, -\pi \, \leq \,
\delta+2d_1\arctan\left(\lambda x\right)-2d_2\arctan\left(x/\lambda\right) \, , $$
once again, forcing $r$ to be very close to $1$ if a small separation is to happen.

A similar analysis can be carried through to obtain information about the separation between consecutive convex or
concave pre-vertices in the case of exterior mappings. The resulting inequalities are analogous to (\ref{sep-convex}) and
(\ref{sep-concave}), with the single term $\delta$ replaced by $2\delta$. The proof will be omitted.

\section{Location of Zeros}

In this section we study the location of the zeros of the Blaschke products appearing in the
representation formulas (\ref{Pf}) and (\ref{zPf}) of Schwarz-Christoffel mappings. Convex or concave mappings
impose no restriction on the location of the zeros, since in the absence of the Blaschke product $B_2$,
(\ref{Pf}) and (\ref{zPf}) will always render univalent mappings. It is probably an ambitious task to determine conditions
on $B_1$ and $B_2$ that are both necessary and sufficient for all mappings of the form (\ref{Pf}) and (\ref{zPf}) to be univalent.
Nevertheless, some necessary conditions can be established. We deal first with the
case of mappings arising from (\ref{Pf}). Because $1+\textrm{Re}\{zf''/f'\}$ will be positive or
negative according to whether $|zB_1|<|B_2|$ or $|zB_1|>|B_2|$, it follows readily from
the radius of convexity for the class $\mathcal{S}$, that we must have
$$ |zB_1(z)|<|B_2(z)| \;\, , \;\, |z|\leq 2-\sqrt{3} \, .$$

\me\noi
In particular, all zeros of $B_2$ must lie in the region $|z|>2-\sqrt{3}$.

\bi\noi{\bf Theorem 9:} {\sl Let $f$ be given by (\ref{Pf}), with $d_1, d_2$ the degrees of the Blaschke products $B_1, B_2$, respectively,
and suppose that $d_2\geq 1$. Suppose that all the zeros of $B_1, B_2$ are contained in the subdisk $|z|\leq r<1$.
Then

\begin{equation}\label{radius} r\,\geq\, \max\left\{\frac{\sqrt{4d_1d_2+9}+3-2d_2}{\sqrt{4d_1d_2+9}+3+2d_2}\, , \,
\frac{2d_2-1-\sqrt{1+4d_1d_2}}{2d_2+1+\sqrt{1+4d_1d_2}}\right\}\, .\end{equation}

\sm
\noi
In particular, if $d_2=1$ then
\begin{equation}\label{uno} r\,\geq\,\frac{\sqrt{4n+5}+1}{\sqrt{4n+5}+5}\,\geq\,\frac12
 \, .\end{equation}
The estimate (\ref{uno}) is sharp for the Koebe function.}

\pf Recall the boundary function $\varphi(t)$ in (\ref{boundary-1}).
At a concave pre-vertex $e^{it_0}$, the exterior angle $2\pi\beta_0\in[-\pi,0)$, and hence $\beta_0\in[-\frac12,0)$.
It follows from (\ref{angle}) that
$$\varphi'(t_0)\leq -2 \, .$$ If we write
$$ B_1(z)=c_1\prod_{k=1}^{d_1}\frac{z-a_k}{1-\overline{a_k}z}  \quad , \quad B_2(z)=c_2\prod_{k=1}^{d_2}\frac{z-b_k}{1-\overline{b_k}z} \, ,$$
then
$$ \varphi'(t)=1+\sum_{k=1}^{d_1}\frac{1-|a_k|^2}{|e^{it}-a_k|^2}-\sum_{k=1}^{d_2}\frac{1-|b_k|^2}{|e^{it}-b_k|^2} \, .$$
After evaluating at $t=t_0$, a simple estimate gives
$$ 1+d_1\frac{1-r}{1+r}-d_2\frac{1+r}{1-r}\leq -2 \, .$$ With $s=(1+r)/(1-r)$ we obtain
$$ d_2s^2-3s-d_1 \geq 0 \, , $$
which implies (\ref{radius}). If $d_2=1$ then $d_1=n-1$, which proves the first estimate in (\ref{uno}).

In order to obtain the second estimate, we observe that at a convex pre-vertex $e^{it_1}$, the exterior angle $2\pi\beta_1>0$,
and therefore $\varphi'(t_1)>0$. This now gives
$$1+d_1\frac{1+r}{1-r}-d_2\frac{1-r}{1+r}>0 \, ,$$
and the second estimate follows.

\sm

To show sharpness, we consider the Koebe function $k(z)=z/(1-z)^2$. Then $k'(z)=(1+z)/(1-z)^3$ and thus
$$ \frac{k''}{k'}(z)=\frac{1}{z+1}-\frac{3}{z-1} \, ,$$
which is consistent with a polygonal mapping onto a 2-gon with a concave vertex with exterior angle $-\pi$ at $k(-1)=-\frac14$,
and a convex vertex with exterior angle $3\pi$ at $k(1)=\infty$ . A calculation gives
$$\frac{k''}{k'}(z)=\frac{1/B_2(z)}{1-z/B_2(z)} \, , $$ with
$B_2(z)=(z+\frac12)/(1+\frac12z)$. Then $r=\frac12$, which coincides with the lower bound in (\ref{uno})
with $n=1$.

\me
\noi
{\bf Remark:} The first estimate in (\ref{radius}) is the better one when $d_1>>d_2$, while the second will
provide better information for $d_2>>d_1$.

\bi

The final theorem describes the analogous situation for mappings of the form (\ref{zPf}) onto the
complement of polygons. The corresponding boundary function is now given by
$$\varphi(t)=\arg\left\{e^{2it}{\frac{B_1}{B_2}(e^{it})}\right\}\, ,$$
for which
$$\varphi'(t)=2+|B_1'(e^{it})|-|B_2'(e^{it})| \, .$$
Since the proof is based on an almost identical analysis, it will be omitted.

\bi\noi{\bf Theorem 10:} {\sl Let $f$ be given by (\ref{zPf}), with $d_1, d_2$ the degrees of the Blaschke products $B_1, B_2$, respectively,
and suppose that $d_2\geq 1$. Suppose that all the zeros of $B_1, B_2$ are contained in the subdisk $|z|\leq r<1$.
Then

\begin{equation}\label{radius-1} r\,\geq\, \max\left\{\frac{\sqrt{d_1d_2+4}+2-d_2}{\sqrt{d_1d_2+4}+2+d_2}\, , \,
\frac{d_2-1-\sqrt{1+d_1d_2}}{d_2+1+\sqrt{1+d_1d_2}}\right\}\, .\end{equation}}

\bi
\sm{\noi}{\bf Acknowledgement:} We thank the referee for a very careful reading of the manuscript
and several valuable suggestions.

\sm

\bibliographystyle{plain}


\me
\noi

\email{}

\end{document}